\documentclass[10pt,oneside,english]{article}
\usepackage[english]{babel}
\usepackage[utf8]{inputenc}
\usepackage{graphicx}
\usepackage{graphics}
\usepackage{geometry}
\usepackage{amssymb}
\usepackage{amsmath, amsthm}
\usepackage{amsmath,amscd}
\usepackage[all,cmtip]{xy}
\usepackage{hyperref}
\usepackage{color}
\usepackage{float}
\usepackage{listings}
\usepackage{subfig}
\usepackage{bm}
\usepackage[final]{pdfpages}
\usepackage{fancyhdr}
\usepackage{epstopdf}
\usepackage{url}
\usepackage{appendix}
\usepackage{verbatim}
\usepackage{color}

\makeatletter
 \theoremstyle{plain}
 
 \theoremstyle{definition}

 \theoremstyle{remark}

\newcommand {\R}{\mathbb{R}}

\renewcommand {\prod} {\stackrel{.}{\wedge}}

\newcommand{\J}{\mathbf{J}}

\newcommand{\diag}{\operatorname{diag}}
\newcommand{\PP}{\mathcal{P}}
\newcommand{\w}{\omega}
\newcommand{\restr}[1]{\vrule height3ex width.4pt depth1.4ex\lower1.4ex\hbox{\scriptsize $\,#1$}}
\newcommand{\rrestr}[1]{\vrule height2ex width.4pt depth0.9ex\lower0.9ex\hbox{\scriptsize $\,#1$}}
\newcommand{\II}{\mathbb{I}}
\newcommand{\ed}{\mathbf{d}}   

\newcommand {\g}{\mathfrak{g}}

\newcommand {\m}{\mathfrak{m}}

\newcommand{\Id}{\text{Id}}

\newcommand{\e}{\mathbf{e}}
\renewcommand{\u}{\mathbf{u}}
\date{}

\begin{document}
\title{The Fast-Superfast Transition in the Sleeping Lagrange Top}
\author{Tudor S. Ratiu\thanks{\texttt{tudor.ratiu@epfl.ch}.
Section de Math\'ematiques, Ecole Polytechnique F\'ed\'erale
de Lausanne,  Lausanne, CH-1015 Switzerland and Mathematical
Sciences Program, Skolkovo Institute of Science and
Technology, 100 Novaya str., Skolkovo, Odintsovsky District,
Moscow Region, Russia 143025.}, Miguel Rodr\'{\i}guez-Olmos
\thanks{\texttt{miguel.rodriguez.olmos@ma4.upc.edu}. Department
of Applied Mathematics IV. Technical University of Catalonia.
Barcelona, Spain.}, Miguel Teixid\'o-Rom\'an
\thanks{\texttt{miguel.teixido@ma4.upc.edu}. Department of Applied
Mathematics IV. Technical University of Catalonia. Barcelona,
Spain.}}

\maketitle
\abstract{The fast-superfast transition is a particular movement
of eigenvalues {found in \cite{lewis1992heavy}} when studying the
family of sleeping equilibria in the Lagrange top. Although
this behaviour of eigenvalues typically suggests a change in
stability or a bifurcation, in this case there is no particular
change in the qualitative dynamical properties of the system.
Using modern methods, based on the singularities of symmetry groups
for Hamiltonian systems, we clarify the appearance of this
transition.}

\
\\

\noindent MSC 2010: 70H05; 70H14; 37J20
\

\section{Introduction}

The Lagrange top is one of the most well-known simple mechanical
systems with symmetries. It consists of an axisymmetric rigid body
with a fixed point moving under the influence of a constant
gravitational force. One of the important characteristics of this
system is that it has  {two symmetry groups:} the
rotations around the axis of gravity and the spin about the axis of
symmetry of the top. In fact, the energy and the conserved quantities
associated with these symmetries make the Lagrange top a Liouville
integrable system (see, e.g.,  \cite{cushman2004global}).

A special configuration of the system is the \emph{sleeping top},
when the top is pointing upwards, so that the axis of gravity and
the symmetry axis coincide. This position has  non-trivial isotropy,
since there is a continuous group of  symmetries that leaves this
configuration invariant. This fact is ultimately responsible for
many of the dynamical properties of a sleeping top. In this work,
we focus on  {this particular configuration, but with
all possible angular velocities. These solutions are
examples of} relative equilibria, motions of
the system which  evolve within the symmetry direction. Generally
speaking, relative equilibria act as organising centres of the
dynamics of a Hamiltonian system with symmetries and, therefore,
their study gives important information about the qualitative
behaviour of the dynamical flow which, typically, is impossible to
obtain analytically.
 
Unlike for generic flows, relative equilibria of Hamiltonian
systems can be characterised as critical points of a certain
function  {depending on a parameter}. This  approach
is very useful, because it gives
necessary conditions for bifurcations of branches of relative
equilibria, based on singularity and critical point theory.
Indeed, when the second variation of this  {parameter dependent} function becomes
degenerate, a new branch of relative equilibria may bifurcate.
Typically, one expects that at the bifurcation point, a transfer
of stability will occur from the original to the bifurcating
branch. However, in the Lagrange top, due to the existence of
continuous isotropy, this is not at all the case. For a large part
of the family of sleeping Lagrange tops, every point is a
bifurcation point to precessing solutions, without the original
branch loosing stability.  The possibility of a bifurcation from a stable branch to a branch that
is also stable was already pointed out by Cartan \cite{Cartan1928}.

\begin{figure}[H]
\begin{center}
\includegraphics[width=13cm]{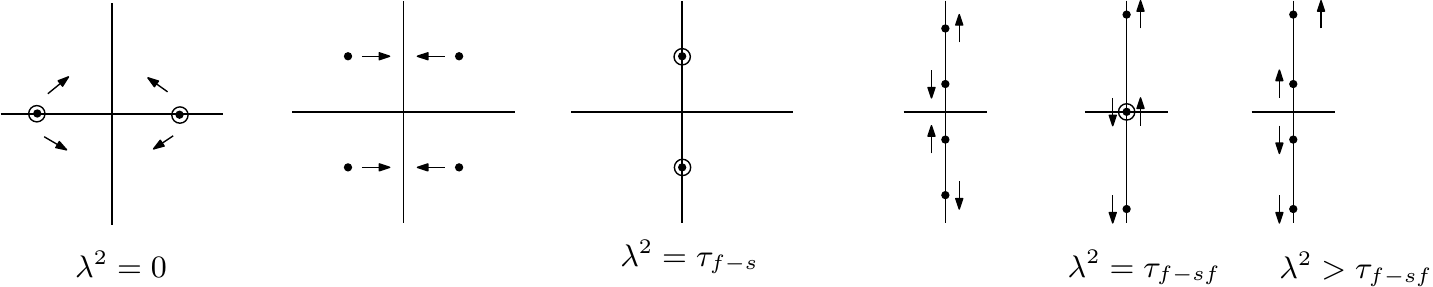}
\caption{Eigenvalues of the sleeping top as a function of the angular velocity $\lambda$ in an oblate top \label{eigen-oblate}}
\end{center}
\end{figure}

These and other phenomena are well known and have been studied in
depth from a geometric perspective in \cite{lewis1992heavy},
 {where the stability range of sleeping Lagrange tops and their
possible bifurcations are computed and the spectral analysis of
the linearised dynamics along the sleeping Lagrange top family
is carried out. Assuming an oblate body, i.e., the largest moment
of inertia is along the symmetry axis, in \cite{lewis1992heavy},
a diagram similar to the one showed in Figure \ref{eigen-oblate}
is presented.}
When the angular velocity $\lambda$ is zero, the eigenvalues are
a real pair. As $\lambda$ is increased, the eigenvalues form a
quadruple and when the angular velocity attains the \emph{fast-slow}
critical value $\tau_{f-s}$, the eigenvalues move onto the
imaginary axis and the system becomes stable. This fast-slow
transition corresponds to a Hamiltonian-Hopf bifurcation (see
\cite{MR1047475}). The surprising fact noticed in
\cite{lewis1992heavy} is that if the angular velocity is further
increased, two of the eigenvalues cross at zero moving along the
imaginary axis for $\lambda^2=\tau_{f-sf}$. This point was called
the \emph{fast-superfast} transition. Usually, such an eigenvalue
crossing implies the existence of a bifurcation or a change in
stability, but from a dynamical point of view, this fast-superfast
point has no special behaviour with respect to nearby points in
the sleeping top family. On the other hand, one could argue that
 since all the points in the stable regime exhibit a bifurcation
to the precessing branch,  a zero eigenvalue at each point in the stable range is expected. However, only a zero
eigenvalue is  observed in the linearisation at the fast-superfast transition point.

\begin{figure}[H]
\begin{center}
\includegraphics[width=13cm]{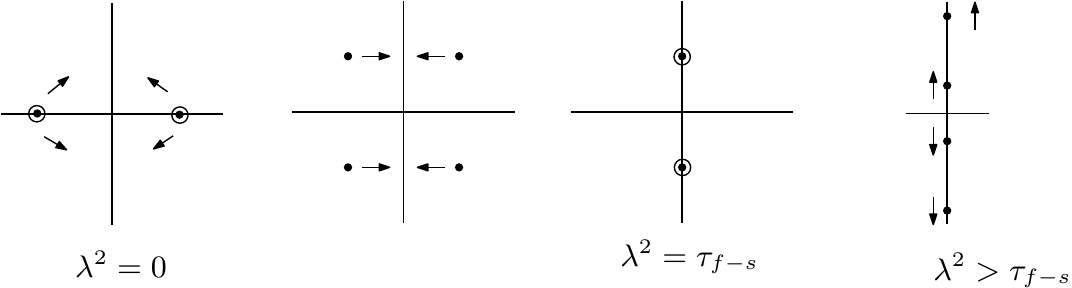}
\caption{Eigenvalues of the sleeping top as a function of the
angular velocity $\lambda$ in prolate top \label{eigen-prolate}}
\end{center}
\end{figure}

If the body is prolate, i.e., the shortest moment of inertia is
along the symmetry axis, the evolution of the eigenvalues is
shown in Figure \ref{eigen-prolate}.
The motions of eigenvalues shown in both figures correspond to those obtained in Figure 3 of \cite{lewis1992heavy}.

As in the oblate case, the fast-slow transition is a Hamiltonian-Hopf
bifurcation between the stable and unstable regimes, but note that
in this case, no fast-superfast transition occurs. It can
be shown that the qualitative dynamical behaviour of the system
does not depend on the oblateness or prolateness of the body, so
 the linearisation should not behave differently
in the oblate and prolate cases. In particular, the existence of
a fast-superfast transition at a unique point in the oblate but
not the prolate case, questions our understanding of the qualitative
dynamics of the sleeping Lagrange top family.

In this paper, we show that this fast-superfast transition is
an artefact due to the existence of continuous isotropy in the
sleeping Lagrange top. This has as a consequence that the
linearisation done in \cite{MR1047475} has implicitly chosen one
angular velocity representative among an infinite number of
possibilities for the equilibria under study. In fact, we show
that for each stable sleeping equilibrium, we can choose the
angular velocity in such a way that the linearisation has a
double zero eigenvalue crossing and therefore exhibits a
fast-superfast transition, for both the oblate and prolate cases.

\section{The Setup}
The Lagrange top is the symmetric Hamiltonian system defined by
the tuple
$$
(T^*\mathrm{SO}(3),\omega_c,\mathbb{T}^2, \J,h),
$$
where the phase space is the cotangent bundle $T^*\mathrm{SO}(3)$
of the proper rotation group $\mathrm{SO}(3)$, equipped with
its canonical symplectic form $\omega_c$. If we use the
identification, given by  the right translation
\begin{eqnarray}
\nonumber
\PP:=\mathrm{SO}(3)\times\R^3 & \rightarrow & T^*\mathrm{SO}(3)\\
\label{rightriv}
(\Lambda,\pi) & \rightarrow &  (\Lambda,\hat\pi\Lambda),\quad
\pi\in\mathbb{R}^3,
\end{eqnarray}
then the toral symmetry group $\mathbb{T}^2$ acts on  on $\PP$ as
 \begin{equation}
 \label{action}  
 ((\theta_1,\theta_2),(\Lambda,\pi))\mapsto
(e^{\theta_1 \widehat{\e_3}}\Lambda e^{-\theta_2 \widehat{\e_3}}, e^{\theta_1\widehat{\e_3}}\pi),\quad (\theta_1,\theta_2)\in
S^1\times S^1.
\end{equation}

 The associated canonical momentum map for this action is given by
\begin{equation}
\label{momentum}
\J(\Lambda,\pi)=(\pi\cdot \e_3,-\pi\cdot \Lambda \e_3)
\in\mathbb{R}^2.
\end{equation}
Finally, the Hamiltonian function $h: \mathcal{P}\rightarrow \mathbb{R}$ is
\[
h(\Lambda,\pi):=mgl \e_3\cdot\Lambda \e_3+
\frac{1}{2}\pi\cdot \mathbb{I}^{-1}_\Lambda \pi,
\]
where $\mathbb{I}_\Lambda:=\Lambda\mathbb{I}\Lambda^{-1}$ and
$\mathbb I:=\diag (I_1,I_1,I_3)$ is the inertia tensor of the
rigid body in the reference configuration with respect to a
principal axes body frame. Notice that, in this reference
configuration, the body symmetry axis is $\e_3$.
Throughout this article we are using both $\mathbb{R}^3$ and
the space $\mathfrak{so}(3)$ of $3\times 3$ antisymmetric matrices
to represent the Lie algebra of $\mathrm{SO}(3)$, identified via
the usual Lie algebra isomorphism given by the hat map
\begin{eqnarray*}
\hat\,:\mathbb{R}^3 & \rightarrow & \mathfrak{so}(3)\\
\mathbf{x}=(x_1,x_2,x_3) & \mapsto & \hat{\mathbf{x}}=
\left[
\begin{array}{ccc}
0 & -x_3 & x_2\\ x_3 & 0 & -x_1\\ -x_2 &  x_1 & 0
\end{array}
\right].
\end{eqnarray*}

In right trivialisation we have $\delta\Lambda=
\widehat{\delta\theta}\Lambda\in T_\Lambda\mathrm{SO}(3)$.
Therefore, a tangent vector $v\in T_{(\Lambda,\pi)}\PP$ is
represented  as $(\delta\theta,\delta\pi)\in
\mathbb{R}^3\times \mathbb{R}^3$.

From this, it is easy to check that the expression for the
symplectic form $\omega$ on $\PP$ (the pull back of $\omega_c$
to $\mathcal{P}$ given by \eqref{rightriv}) is given by
\begin{equation}
\label{symp}
\w(\Lambda,\pi)((\delta\theta_1,\delta\pi_1)
(\delta\theta_1,\delta\pi_2))=
\delta\pi_2\cdot\delta\theta_1-\delta\pi_1\cdot\delta\theta_2-
\pi\cdot(\delta\theta_1\times\delta\theta_2).
\end{equation}

\section{The sleeping Lagrange top}
\label{sec:sleeping_Lagrange}
In this model, the sleeping Lagrange top is a relative equilibrium
at the phase space point $z=(\Lambda,\pi):=
(\mathrm{Id},\lambda I_3\e_3)$ with arbitrary velocity
$(\xi_L,\xi_R)\in\mathbb{R}^2=\mathrm{Lie}\,\mathbb{T}^2$
satisfying $\xi_L-\xi_R=\lambda$. This makes the set of all
admissible velocities at each relative equilibrium a
one-parameter family. In order to see this, and using standard
arguments from geometric mechanics (see, e.g., \cite{Mars92},
\cite{MaRa1999}), we have to check that the augmented Hamiltonian
\begin{equation}
\label{augham}
h_{(\xi_L,\xi_R)}(\Lambda,\pi):=mgl \e_3\cdot\Lambda \e_3+
\frac{1}{2}\pi\cdot \mathbb{I}^{-1}_\Lambda \pi-
\pi\cdot\e_3\xi_L+\pi\cdot\Lambda\e_3\xi_R
\end{equation}
has a critical point at $(\mathrm{Id},\lambda I_3\e_3)$. The
derivative of $h_{(\xi_L,\xi_R)}$ at an arbitrary point
$(\Lambda,\pi) \in \mathcal{P}$ is
\[
\ed h_{(\xi_L,\xi_R)} (\Lambda,\pi)(\delta\theta,\delta\pi)=
mgl\e_3\cdot (\delta\theta\times \Lambda\e_3)+
\delta\pi\cdot \II_\Lambda^{-1} \pi+
\pi\cdot (\delta\theta\times\II_\Lambda^{-1} \pi)-
\xi_L\delta\pi\cdot \e_3+\xi_R\delta\pi\cdot \Lambda\e_3+
\xi_R\pi\cdot(\delta\theta\times\Lambda\e_3)
\]
At $(\Lambda,\pi)=(\mathrm{Id},\lambda I_3\e_3)$ this becomes
\begin{equation}
\label{1var}
\ed h_{(\xi_L,\xi_R)}(\mathrm{Id},\lambda I_3\e_3)
(\delta\theta,\delta\pi)=
(\lambda-(\xi_L-\xi_R))\delta \pi\cdot \e_3
\end{equation}
which vanishes precisely when $\xi_L-\xi_R=\lambda$. Therefore,
according to \eqref{action},  the dynamical evolution of
$z=(\mathrm{Id},\lambda I_3\e_3)$ is given by $z(t)=
(e^{t\lambda\widehat{\e_3}},\lambda I_3\e_3)$.

In order to study the linearisation of the Hamiltonian system at
the sleeping equilibrium, we need to compute several more geometric
objects. First, from \eqref{momentum}, it is clear that
$$
\mu:=\J(z)= (\lambda I_3,-\lambda I_3).
$$
 Second, since $\mathbb{T}^2$ is Abelian, the
coadjoint stabiliser of $\mu$ is $G_\mu=\mathbb{T}^2$. Third,
from \eqref{action}, we find that the stabiliser of the phase
space point $z=(\mathrm{Id},\lambda I_3\e_3)$ is
$$
G_z = {S^1}^D:=\{(\theta,\theta)\in\mathbb{T}^2\,:\,
\theta\in S^1\}.
$$
A normalised basis for its Lie algebra $\mathfrak{g}_z$
is $\frac{1}{\sqrt{2}}(1,1)$, i.e.,
$$
\g_z  =  \mathrm{span}\left\{\frac{1}{\sqrt{2}}(1,1)\right\} .
$$

Define $\m:=\mathrm{span}\left\{\frac{1}{\sqrt{2}}(1,-1)\right\}$ which is a complement to $\g_z$ in $\g$, i.e.,
\begin{equation}
\label{gmusplit}
\g=\m\oplus\g_z.
\end{equation}
According to this direct sum decompostion, the velocity of the
sleeping Lagrange top takes the form
\begin{equation}
\label{ortvel}
(\xi_L,\xi_R)=\frac\lambda2(1,-1)+\eta(1,1),
\end{equation}
where $\eta \in \mathbb{R}$ is arbitrary.  Since $G_z$ is a
continuous subgroup of positive dimension, \eqref{ortvel}
reflects the fact that the velocity of the relative equilibrium
$z$, for which the first variation \eqref{1var} vanishes, is
defined only up to an element of $\g_z$. However, its projection
$\xi^\perp$ onto the subspace $\m$, according to the splitting
\eqref{gmusplit}, sometimes called the \emph{orthogonal velocity}
of the equilibrium, is unique.

\section{The Symplectic Slice}
The linearisation of the Hamiltonian system at the relative
equilibrium $z=(\mathrm{Id},\lambda I_3\e_3)$ is given by
\begin{equation}
\label{linear}
L:=\Omega_N^{-1}\ed^2h_{(\xi_L,\xi_R)}
(\mathrm{Id},\lambda I_3\e_3)\restr{N}\in
\mathfrak{sp}(N),
\end{equation}
where $(N,\Omega_N)$ is a symplectic vector subspace of $T_z\PP$
defined by an arbitrary $G_z$-invariant splitting
$$
\ker T_z\J = T_z(G_\mu\cdot z)\oplus N,
$$
$\Omega_N=\w(z)\rrestr{N}$, and $\ed^2h_{(\xi_L,\xi_R)}
(\mathrm{Id},\lambda I_3\e_3)$ denotes the Hessian of
$h_{(\xi_L,\xi_R)}$ at the equilibrium point
$(\mathrm{Id},\lambda I_3\e_3)$. The vector space
$(N,\Omega_N)$ is often called a \emph{symplectic slice}.

In order to compute the symplectic slice, we start by studying
the fundamental vector fields $(\xi_1,\xi_2)_{\PP}$ for the
action \eqref{action}, which are given, in  the right
trivialisation \eqref{rightriv}, by
$$
(\xi_1,\xi_2)_{\PP}(\Lambda,\pi)=(\xi_1\e_3-\xi_2\Lambda\e_3,\xi_1\e_3\times\pi).
$$
So, at the sleeping equilibrium $(\Id,\lambda I_3\e_3)$, we have
$$
T_{(\Id,\lambda I_3\e_3)}(G_\mu\cdot z) = \mathrm{span}\left\{(\e_3,0)\right\},
$$
where $\mu=\J(z)= (\lambda I_3,-\lambda I_3)$ and the coajoint
isotropy subgroup is $G_ \mu = \mathbb{T}^2$ (see
\S\ref{sec:sleeping_Lagrange}).

Since the derivative of the momentum map is
$$
T_{(\Lambda,\pi)}\J(\delta\theta,\delta\pi)=
(\delta\pi\cdot\e_3,-\delta\pi\cdot\Lambda\e_3-
\pi\cdot(\delta\theta\times\Lambda\e_3)),
$$
at the sleeping Lagrange top point $({\rm Id},
\lambda I_3 \mathbf{e}_3)$ this becomes
$$
T_{(\Id,\lambda I_3\e_3)}\J(\delta\theta,\delta\pi)=
(\delta\pi\cdot\e_3,-\delta\pi\cdot\e_3-
\lambda I_3\e_3\cdot(\delta\theta\times\e_3))
$$
and, therefore,
$$
\ker T_{(\Id,\lambda I_3\e_3)}\J=\{(\delta\theta,\delta\pi)\,:\,\delta\pi\perp\e_3\}.
$$
Hence, a possible choice for the symplectic slice at $z$ is
$$N=\{(\delta\theta,\delta\pi)\,:\,\delta\theta,\delta\pi\perp\e_3\},$$
with symplectic form
\begin{equation}\label{Lsymp}\Omega_N\left((\delta\theta_1,\delta\pi_1),(\delta\theta_2,\delta\pi_2)   \right)=\delta\pi_2\cdot\delta\theta_1-\delta\pi_1\cdot\delta\theta_2-\lambda I_3\e_3\cdot(\delta\theta_1\times\delta\theta_2).\end{equation}
It can easily be checked that $N$ is indeed $G_z$-invariant and
that  the isotropy group $G_z={S^1}^D$ acts on $N$ by
$$
\varphi\cdot(\delta\theta,\delta\pi)=
(e^{{\rm i}\varphi\widehat{\e_3}}\delta\theta,
e^{{\rm i}\varphi\widehat{\e_3}}\delta\pi), \qquad
\varphi \in S^1 \cong {S^1}^D.
$$

\section{Stability}
In this section, we briefly reproduce the classic stability
result for a sleeping equilibrium of the Lagrange top in a
formulation adapted to the general framework of this article.
This is necessary for the subsequent bifurcation analysis relevant
to the fast-superfast transition.

The  second variation of the augmented Hamiltonian \eqref{augham} is
\begin{eqnarray*}
\ed^2 h(\xi_L,\xi_R) ((\delta\theta_1,\delta\pi_1),
(\delta\theta_2,\delta\pi_2)) & = &
(-mgl-\lambda^2I_3+\lambda^2\frac{I_3^2}{I_1}-
\lambda I_3\xi_R)(\delta\theta_1\times\e_3)\cdot
(\delta\theta_2\times\e_3)\\
& + & \left(\lambda\left(1-\frac{I_3}{I_1}\right)+\xi_R\right)\left[\delta\pi_1\cdot(\delta\theta_2\times\e_3)+
\delta\pi_2\cdot(\delta\theta_1\times\e_3)\right]\\
& + & \delta\pi_1\cdot\II^{-1}\delta\pi_2.
\end{eqnarray*}

Let
\begin{equation}
\label{basis}
\u_1=(\e_1,0),\,\u_2=(\e_2,0),\,\u_3=(0,\e_1),\,\u_4=(0,\e_2)
\end{equation}
be a basis of $N$. Using it and \eqref{ortvel}, we get $\xi_R=\eta-\lambda/2$. Thus, at the sleeping top equilibrium,
\begin{equation}
\label{stabform}
\ed^2 h_{(\xi_L,\xi_R)}(\Id,\lambda I_3\e_3)\restr{N}
=\begin{bmatrix}
A & 0 & 0 & B \\
0 & A & -B & 0 \\
0 & -B & C & 0 \\
B & 0 & 0 & C
\end{bmatrix} ,
\end{equation}
where
\[
A:=-mgl+\frac{\lambda^2I_3}{2I_1}(2I_3-I_1)-\lambda I_3\eta,\quad
B:=\frac{\lambda}{2I_1}(2I_3-I_1)-\eta,\quad
C:=\frac{1}{I_1}.
\]
The eigenvalues of this matrix are
\[
\sigma_{\pm}:=\frac{1}{2}\left((A+C)\pm\sqrt{(A+C)^2-4(AC-B^2)}
\right).
\]
Thus, $\sigma_{\pm}>0\Longleftrightarrow AC-B^2>0$.

Recall from the theory of stability of Hamiltonian relative
equilibria (see \cite{MonRO11} for a reference appropriate to
the isotropy-based approach consistent with this article),
that in order to guarantee nonlinear stability it is enough
to find an admissible velocity such that \eqref{stabform} is
definite.  For a sleeping top equilibrium of the form
$z=(\Id,\lambda I_3\e_3)$, this is equivalent to saying that
it is stable if we can find $\eta\in\R$ such that $AC-B^2>0$.

Note that
\[
AC-B^2=-\eta^2+\frac{I_3-I_1}{I_1}\lambda\eta+
\frac{I_3\lambda^2}{2I_1}-\frac{mgl}{I_1}-\frac{\lambda^2}{4}.
\]
Therefore, for a fixed value of $\lambda$, this is a quadratic
function of $\eta$ that attains its maximum value when $\eta$
is equal to
\[
\eta^*:=\frac{I_3-I_1}{2I_1}\lambda.
\]
Substituting this value of $\eta$ in the previous expression yields
\[
AC-B^2\restr{\eta=\eta^*}=
\frac{\lambda^2I_3^2-4mglI_1}{4I_1^2}
\]
and we  conclude that if
\[\lambda^2> \frac{4mglI_1}{I_3^2}\]
then the sleeping top equilibrium $(\mathrm{Id},
\lambda I_3 \mathbf{e}_3)\in \PP$ is nonlinearly stable. This
is the classical fast top condition. Notice that if
\[
\lambda^2< \frac{4mglI_1}{I_3^2}
\]
then the stability test is inconclusive and this regime must be
studied by linearisation methods, as we do in the next section.

\section{Linearisation}
Formula \eqref{Lsymp} easily implies that the expressions of the
linear symplectic form  $\Omega_N$ and its inverse, relative to
the basis \eqref{basis},  are
\[
\Omega_N=\begin{bmatrix}
0 & -\lambda I_3 & 1 & 0 \\
\lambda I_3 & 0 & 0 & 1 \\
-1 & 0 & 0 & 0 \\
0 & -1 &  & 0
\end{bmatrix},\qquad
\Omega_N^{-1} =
\begin{bmatrix}
0 & 0 & -1 & 0\\
0 & 0 & 0 & -1\\
1 & 0 & 0 & -\lambda I_3 \\
0 & 1 & \lambda I_3 & 0
\end{bmatrix}
\]
so the matrix of the linearised system \eqref{linear} is
\begin{equation}
\label{Leta}
L=\begin{bmatrix}
     0 & B & -C & 0 \\
     -B & 0 & 0 & -C \\
     A-\lambda I_3B & 0 &  0 & B-\lambda I_3C\\
     0 & A-\lambda I_3B  & -B+\lambda I_3C & 0
    \end{bmatrix}.
\end{equation}
The characteristic polynomial of $L$ is
\[
t^4+\left(2(AC-B^2)+(2B-CI_3\lambda)^2 \right) t^2+(AC-B^2)^2,
\]
which, after some manipulations, can be written as
$$
t^4+\frac{E^2+F}{2I_1^2} t^2+\frac{(E^2-F)^2}{16I_1^4},
$$
where
\[
E := I_3 \lambda - I_1 (2 \eta + \lambda),\qquad
F := I_3^2 \lambda^2-4mg l  I_1.
\]
Thus, the eigenvalues of $L$ are
\begin{equation}
\varrho_{\pm,\pm}:=\frac{i}{2I_1}(\pm E \pm \sqrt{F}).
\label{eigenvalues}
\end{equation}

Notice that for any $\lambda$, that is, for any sleeping top
equilibrium, an admissible velocity (equivalently, a value of
$\eta$) can be chosen so that $E$ takes any real value. Also,
$F$ is positive if and only if $\lambda$ corresponds to an
equilibrium for which there is an $\eta$  making  the matrix
\eqref{stabform} definite. Therefore, we conclude that:
\begin{itemize}
\item For any $\lambda^2<\frac{4mglI_1}{I_3^2}$, the eigenvalues
of the linearisation have non-zero real part and therefore
the sleeping equilibrium is unstable. $E$ can be chosen so that
the imaginary part is zero (real double pair) or non-zero
(complex quadruple) at each point of the unstable regime.
\item For any $\lambda^2>\frac{4mglI_1}{I_3^2}$, the eigenvalues
of the linearisation are purely imaginary. At each point of the
stable range, $E$ can be chosen  so that we have 4 distinct
non-zero eigenvalues, two double imaginary non-zero eigenvalues,
or an imaginary pair and a double zero.
\end{itemize}

Notice that the last property guarantees the existence of one
imaginary pair plus a doble zero of the linearised system at
each point of the stable range, provided a suitable admissible
velocity at that point is chosen. This shows that the
fast-superfast transition actually happens at each point of
the stable regime if we use the freedom in the isotropy Lie
algebra when linearising the system at a sleeping equilibrium.

\medskip

We now study the relationship between our approach and the results
in \cite{lewis1992heavy}. First, we notice that in order to obtain
a zero crossing of eigenvalues along the sleeping top equilibrium
family (the fast-superfast transition) the condition to be
satisfied is
$$
E^2=F,
$$
which is equivalent to
\begin{equation}
\label{hyperbola}
(2I_3-I_1)\lambda^2+4(I_3-I_1)\eta\lambda-4I_1\eta^2-4mgl=0.
\end{equation}
Notice that the moments of inertia inequalities $I_i<I_j+I_k$
 for any distinct $i, j, k$,
imply that this quadric is a hyperbola in the
$(\lambda,\eta)$-plane for any choice of inertia tensor.
As we have seen before, one can choose $\eta$ such that
$(\lambda,\eta)$ lies on the hyperbola given by \eqref{hyperbola}
if and only if $F>0$ (stable sleeping equilibrium), which holds
for both the prolate and oblate cases, as opposed to was observed
in  \cite{lewis1992heavy}.

A reason for this apparent disagreement is the following. The
approach taken in \cite{lewis1992heavy}, and based on reduction of
the right $S^1$-action, corresponds, in our setup, to making the
permanent choice $\eta=\lambda/2$ along the family of sleeping equilibria. Therefore, in \cite{lewis1992heavy}, along this
family,
$$
(\xi_L,\xi_R)=(\lambda,0).
$$
Using this choice of velocities, we conclude that \eqref{hyperbola}
is equivalent to
\begin{equation}\label{fsfcondition}
\lambda^2=\frac{mgl}{I_3-I_1},\end{equation}
which corresponds, in the axisymmetric case $I_1=I_2$, to the
fast-superfast conditions $(4.24)$ of \cite{lewis1992heavy}.
 Notice that, in the $(\lambda,\eta)$-plane, for a prolate top
$I_3<I_1$, the line $\eta=\lambda/2$ never meets the hyperbola
\eqref{hyperbola}. However, for an oblate top $I_3>I_1$, the
fast-superfast transition is found at the unique point
\eqref{fsfcondition}, which is the intersection of this line
with the hyperbola. Both observations are in total agreement
with \cite{lewis1992heavy}.

These considerations are illustrated in the following figures.
In Figure \ref{Foblate}, the different possibilities for the
eigenvalues of the linearised system are shown in the oblate case.
The line $L_1$,  corresponding to the choice $\eta=\lambda/2$ of \cite{lewis1992heavy},  intersects the hyperbola $E^2=F$ precisely
at the value $\lambda=\tau_{f-sf}$. Notice, however, that it is possible to find a different relationship between $\lambda$ and
$\eta$ such that this intersection happens at any value of
$\lambda$ greater than $\tau_{f-s}$, that is, the fast-superfast
transition happens at each point of the stable regime in the family
of sleeping equilibria for oblate tops.

In Figure \ref{Fprolate} the analogous situation for prolate tops
is presented. Notice that in this case, the line $L_1$ corresponding
to the velocity choice in \cite{lewis1992heavy} does not intersect
the hyperbola $E^2=F$. However, as is apparenr from the figure, and
exactly as for the oblate case, a different linear relationship
between $\lambda$ and $\eta$ can be chosen so that this intersection
exists for any value of $\lambda$ in the stable regime.

In both figures we have shown the linear relationship between
$\lambda$ and $\eta$ corresponding to the condition $E=0$.
Along that line, the eigenvalues of the linearisation always
come in pairs; the fast-superfast and the fast-slow transitions
occur at the same point.
\begin{figure}
\centering
\begin{minipage}{.5\textwidth}
  \centering
\includegraphics[width=5cm]{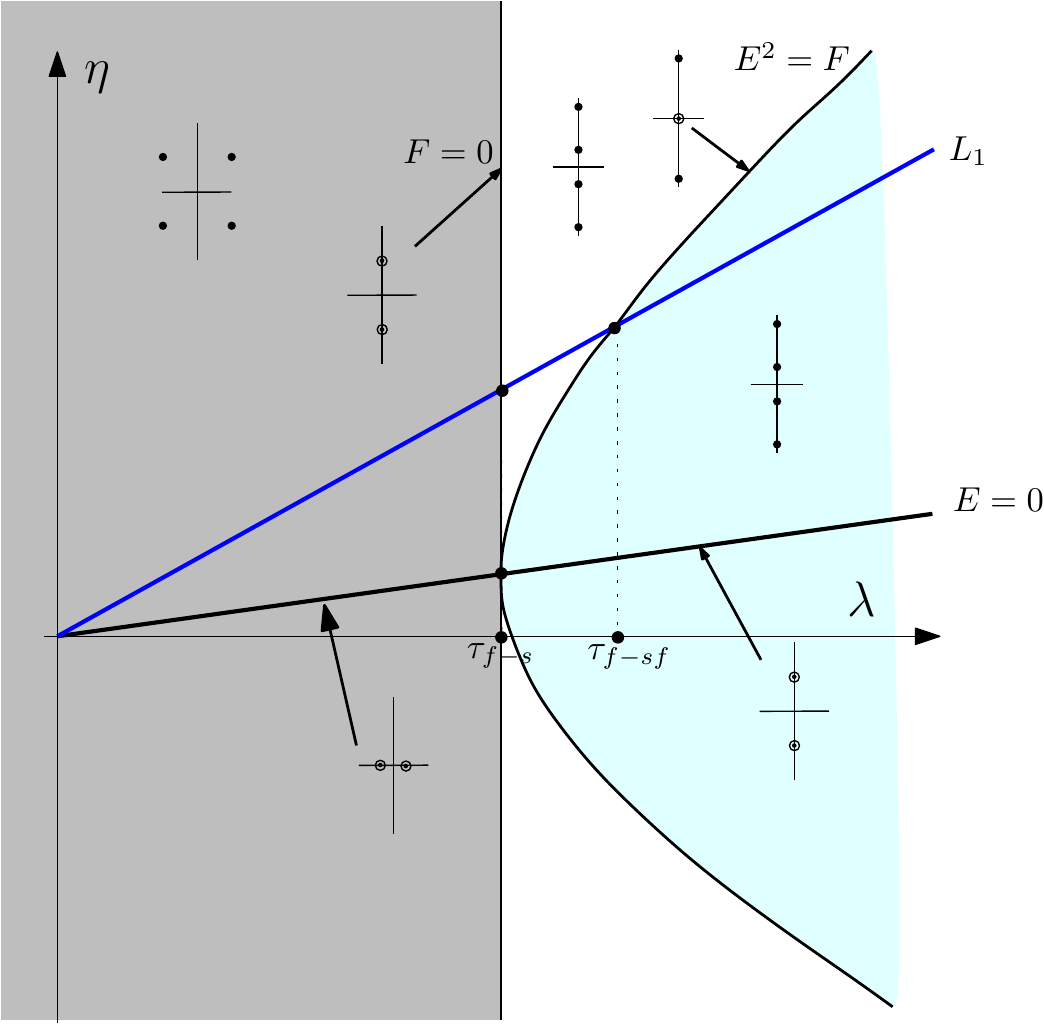}
\caption{Oblate case}\label{Foblate}
\end{minipage}%
\begin{minipage}{.5\textwidth}
  \centering
\includegraphics[width=5cm]{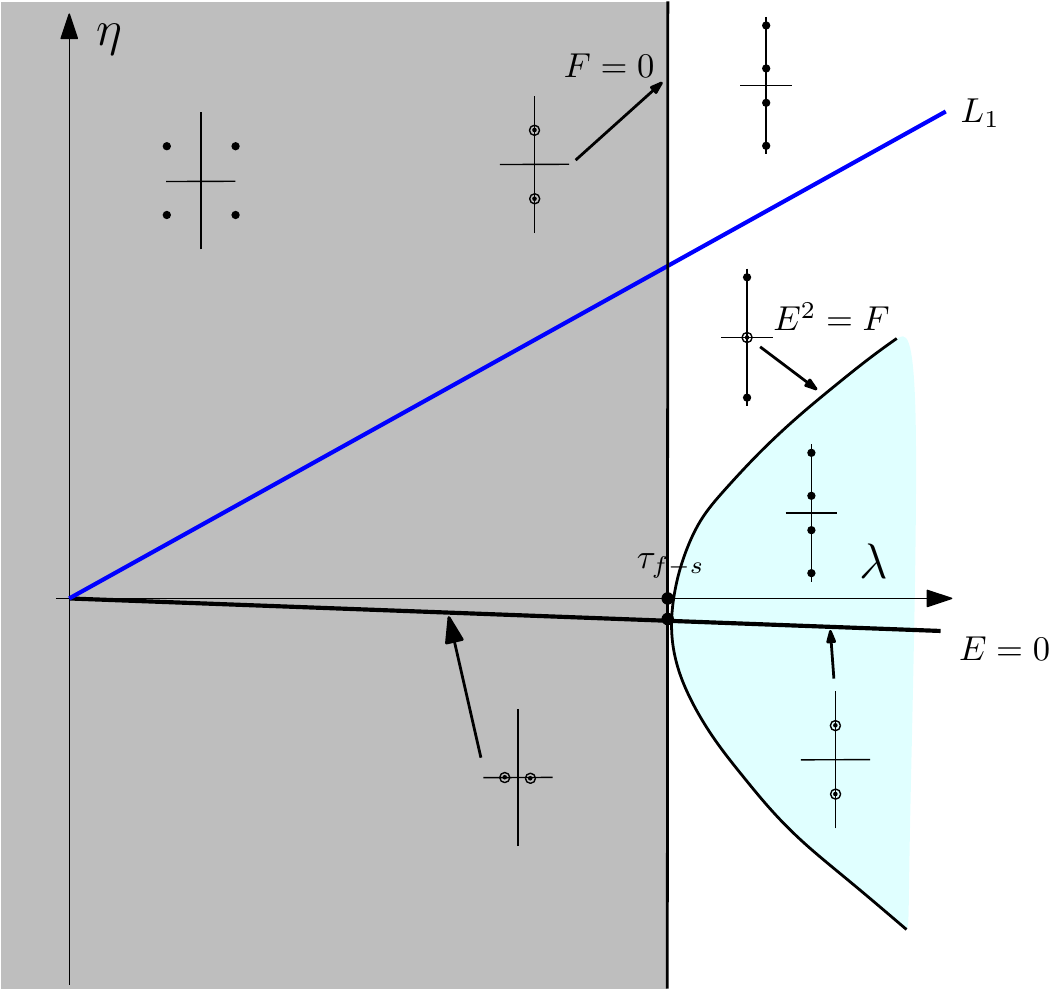}
\caption{Prolate case}\label{Fprolate}
\end{minipage}
\end{figure}

\section{Conclusion}
The existence of continuous isotropy, and therefore an isotropy Lie
algebra $\g_z$ of positive dimension, at a relative equilibrium
reflects the fact that the linearisation of the Hamiltonian vector
field at the equilibrium, given by \eqref{linear}, is not uniquely
defined. The effect of this is that the set of all linearised
vector fields should be considered in order to avoid missing
pieces of information about the qualitative properties of the
Hamiltonian flow near the equilibrium point. In the case of the
sleeping Lagrange top, we have shown that the systematic study
of the $\g_z$-parametrised family of linearisations given by
\eqref{Leta} guarantees that the fast-superfast transition occurs
for both prolate and oblate Lagrange tops. Furthermore, in both
cases, this transition can be observed for all points in the stable
range along the family of sleeping equilibria.

It can be shown that the fast-superfast transitions are closely
related to the bifurcations from sleeping to precessing relative
equilibria, which are known to happen precisely at each point of
the stable range of sleeping tops. This will be studied in detail
in \cite{MonROpp}.

\paragraph{Acknowledgements.}
M. R.-O. and M. T.-R. acknowledge the financial support of the Ministerio de Ciencia e Innovaci\'on (Spain), project
MTM2011-22585 and AGAUR, project 2009 SGR:1338. M. T.-R. thanks
for the support of a FI-Agaur Ph.D. Fellowship. M. R.-O. was
partially supported by the EU-ERG grant ``SILGA". T. S. R. was
partially supported by NCCR SWISSMap and grant 200021-140238,
both of the Swiss National Science Foundation.

{
\end{document}